\documentclass[letterpaper, 10 pt, conference]{ieeeconf}  

\IEEEoverridecommandlockouts                              

\title{\LARGE \bf
The Variational Attitude Estimator in the Presence of Bias in Angular Velocity Measurements
}

\author{Maziar Izadi$^{1}$, Sasi Prabhakaran$^{2}$, Amit Sanyal$^{2,\dag}$, Carlos Silvestre$^{3}$, and Paulo Oliveira$^{4}$
\thanks{$^{1}$M. Izadi is with the Department of Aerospace Engineering, Texas A\&M University, College Station, TX 77840.
        {\tt\small maziar@tamu.edu}}%
\thanks{$^{2}$S.P. Viswanathan and A.K. Sanyal are with the Department of Mechanical \& Aerospace Engineering, Syracuse University, 
Syracuse, NY 13244.
        {\tt\small \{sviswana,aksanyal\}@syr.edu}}%
\thanks{$^{3}$C. Silvestre is with Faculty of Science and Technology, University of Macau, China
        {\tt\small csilvestre@umac.mo}}%
\thanks{$^{4}$P. Oliveira is with ISR - LARSyS - IST - ULisboa, Lisbon, Portugal
        {\tt\small pjcro@isr.ist.utl.pt}}%
\thanks{$^{\dag}$ Address all correspondence to this author.}
}

\usepackage{graphicx}          
\usepackage{amsmath,amsfonts,amssymb}
\usepackage{mathtools, textcomp}
\usepackage{mathrsfs}
\usepackage{graphicx}
\usepackage[dvips]{epsfig}
\usepackage{epstopdf}
\usepackage{color}
\usepackage{dsfont}

\newcommand{\SO}{\ensuremath{\mathrm{SO(3)}}}

\newcommand{\Ta}{\ensuremath{\mathrm{T}}}
\newcommand{\T}{^{\mbox{\small T}}}
\newcommand{\so}{\ensuremath{\mathfrak{so}(3)}}
\newcommand{\SE}{\ensuremath{\mathrm{SE(3)}}}

\newcommand{\bR}{\ensuremath{\mathbb{R}}}
\newcommand{\bS}{\ensuremath{\mathbb{S}}}

\newcommand{\mrm}{\mathrm}

\newcommand{\diag}{\mbox{diag}}
\newcommand{\bbm}{\begin{bmatrix}}
\newcommand{\ebm}{\end{bmatrix}}
\newcommand{\matl}{\left[ \begin{array}}
\newcommand{\matr}{\end{array} \right]}
\newcommand{\be}{\begin{equation}}
\newcommand{\ee}{\end{equation}}
\newcommand{\bea}{\begin{eqnarray}}
\newcommand{\eea}{\end{eqnarray}}
\newcommand{\beas}{\begin{eqnarray*}}
\newcommand{\eeas}{\end{eqnarray*}}
\newcommand{\nn}{\nonumber}

\newcommand{\cL}{\mathcal{L}}

\newcommand{\cT}{\mathcal{T}}
\newcommand{\cU}{\mathcal{U}}

\newcommand{\di}{\mathrm{d}}

\newcommand{\lan}{\langle}
\newcommand{\ran}{\rangle}

\newcommand{\U}{U^{m}}
\newcommand{\cS}{\mathcal{S}}

\DeclareMathOperator{\expm}{{expm}}  			

\newtheorem{theorem}{Theorem}[section]

\newtheorem{proposition}[theorem]{Proposition}

\def\thesection{\arabic{section}}

\newcommand{\bi}{\begin{itemize}}
\newcommand{\ei}{\end{itemize}}

\DeclareMathAlphabet{\mathpzc}{OT1}{pzc}{m}{it}

\begin{document}

\maketitle
\thispagestyle{empty}
\pagestyle{empty}

\begin{abstract}
Estimation of rigid body attitude motion is a long-standing problem of interest in several 
applications. This problem is challenging primarily because rigid body motion is described by 
nonlinear dynamics and the state space is nonlinear. The extended Kalman filter and its several 
variants have remained the standard and most commonly used schemes for attitude estimation 
over the last several decades. These schemes are obtained as approximate solutions to the 
nonlinear optimal filtering problem. However, these approximate or near optimal solutions may 
not give stable estimation schemes in general. The variational attitude estimator was introduced 
recently to fill this gap in stable estimation of arbitrary rigid body attitude motion in the presence of 
uncertainties in initial state and unknown measurement noise. This estimator is obtained by applying 
the Lagrange-d'Alembert principle of variational mechanics to a Lagrangian constructed 
from residuals between measurements and state estimates with a dissipation term that is 
linear in the angular velocity measurement residual. In this work, the variational attitude estimator 
is generalized to include angular velocity measurements that have a constant bias in addition to 
measurement noise. The state estimates converge to true states almost globally over the state 
space. Further, the bias estimates converge to the true bias once the state estimates converge to 
the true states. 
\end{abstract}


\section{Introduction}\label{Sec1}
Estimation of attitude motion is essential in applications to spacecraft, unmanned aerial and 
underwater vehicles as well as formations and networks of such vehicles. In this work, we consider 
estimation of attitude motion of a rigid body from measurements of known inertial directions and 
angular velocity measurements with a constant bias, where all measurements are made with 
body-fixed sensors corrupted by sensor noise. The number of direction vectors measured by the 
body may vary over time. For the theoretical developments in this paper, it is assumed that at 
least two directions are measured at any given instant; this assumption ensures that the attitude 
can be uniquely determined from the measured directions at every instant. The attitude estimation 
scheme presented here follows the variational framework of the estimation scheme recently 
reported in \cite{Automatica,ICRA2015}. Like the estimation scheme in~\cite{Automatica}, the scheme 
presented here has the following important properties: (1) attitude is represented globally over 
the configuration space of rigid body attitude motion without using local coordinates or 
quaternions; (2) no assumption is made on the statistics of the measurement noise; (3) unlike 
model-based estimation schemes (e.g., 
\cite{leishman2014quadrotors,bras2013nonlinear,morgado2014embedded}), no knowledge of 
the attitude dynamics model is assumed; (4) the estimation scheme is obtained by applying the 
Lagrange-d'Alembert principle from variational mechanics~\cite{gold,green} to a Lagrangian 
constructed from the measurement residuals with a dissipation term linear in the angular velocity 
measurement residual; and (5) the estimation scheme is discretized for computer implementation 
by applying the discrete Lagrange-d'Alembert principle~\cite{marswest, haluwa}. It is assumed that measurements of direction vectors and angular velocity 
are available at sufficient frequency, such that a dynamics model is not needed to propagate 
state estimates between measurements.

The earliest solution to the attitude determination problem from two inertial vector measurements 
is the so-called ``TRIAD algorithm" from the early 1960s~\cite{black64}. This was followed by 
developments in the problem of attitude determination from a set of vector measurements, 
which was set up as an optimization problem called Wahba's problem~\cite{jo:wahba}. This 
problem of instantaneous attitude determination has many different solutions  
in the prior literature, a sample of which can be obtained in 
\cite{jo:solwahba,markley1988attitude,sanyal2006optimal}. Much of the published 
literature on estimation of attitude states use local coordinates or unit quaternions
to represent attitude. Local coordinate representations, including commonly used quaternion-derived parameters like
the Rodrigues parameters and the modified Rodrigues parameters (MRPs), cannot describe 
arbitrary or tumbling attitude motion, while the unit quaternion representation of attitude is 
known to be ambiguous. Each physical attitude corresponds to an element of the Lie group 
of rigid body rotations $\SO$, and can be represented by a pair of antipodal quaternions on the 
hypersphere $\bS^3$, which is often represented as an embedded submanifold of $\bR^4$ in 
attitude estimation. For dynamic attitude estimation, this ambiguity in the representation could 
lead to instability of continuous state estimation schemes due to unwinding, as is described 
in~\cite{bhat,CSMpaper,jgcd12}. 

Attitude observers and filtering schemes on $\SO$ and $\SE$ have been reported in, e.g.,  
\cite{sanyal2006optimal,silvest08,Lageman,markSO3,mahapf08,bonmaro09,Vas1}. 
These estimators do not suffer from kinematic singularities like estimators using coordinate 
descriptions of attitude, and they do not suffer from the unstable unwinding phenomenon which 
may be encountered by estimators using unit quaternions. Many of these schemes are based on 
near optimal filtering and do not have provable stability. Related to Kalman filtering-type schemes  
is the maximum-likelihood (minimum energy) filtering method of Mortensen~\cite{Mortensen}, 
which was recently applied to attitude estimation, resulting in a nonlinear attitude estimation 
scheme that seeks to minimize the stored ``energy" in measurement 
errors~\cite{aguhes06,ZamPhD}. 
This scheme is obtained by applying Hamilton-Jacobi-Bellman (HJB) theory~\cite{kirk} to 
the state space of attitude motion, as shown in \cite{ZamPhD}. Since the 
HJB equation can be only approximately solved with increasingly unwieldy expressions for 
higher order approximations, the resulting filter is only ``near optimal" up to second order. 
Unlike the filtering schemes that are based on Kalman filtering or ``near optimal" solutions of the HJB 
equation and do not have provable stability, the estimation scheme obtained here is shown to be almost globally asymptotically stable even in the case of biased 
angular velocity measurements. The special case of unbiased velocity measurements was 
dealt with in a prior version of this estimator that appeared recently~\cite{Automatica}. Moreover, 
unlike filters based on Kalman filtering, the estimator proposed here does not require any 
knowledge about the statistics of the initial state estimate or the sensor noise.   

This paper is structured as follows. Section \ref{Sec3} details the measurement model for measurements of 
inertially-known vectors and biased angular velocity measurements using body-fixed sensors. The problem of variational attitude estimation from these 
measurements in the presence of rate gyro bias is formulated and solved on $\SO$ in Section 
\ref{Sec4}. A Lyapunov stability proof of this estimator is given in Section \ref{Sec5}, along with 
a proof of the almost global domain of convergence of the estimates in the case of perfect 
measurements. It is also shown that the bias estimate converges to the true bias in this case. 
This continuous estimation scheme is discretized in Section \ref{Sec6} in the form of a Lie group 
variational integrator (LGVI) using the discrete 
Lagrange-d'Alembert principle. Numerical simulations are carried out using this LGVI as the discrete-time variational 
attitude estimator in Section \ref{Sec6} with a fixed set of gains. Section \ref{Sec7} gives concluding remarks, contributions and possible future 
extensions of the work presented in this paper.

\section{Measurement Model}\label{Sec3}

For rigid body attitude estimation, assume that some inertially-fixed vectors are measured 
in a body-fixed frame, along with body angular velocity measurements having a constant bias.
Let $k\in\mathbb{N}$ known inertial vectors be measured in a coordinate frame fixed to the 
rigid body. Denote these measured vectors as $u_j^m$ for $j=1,2,\ldots,k$, in the body-fixed 
frame. Denote the corresponding known vectors represented in inertial 
frame as $e_j$; therefore $u_j=R \T e_j$, where $R$ is 
the rotation matrix from the body frame to the inertial frame. This rotation matrix provides a 
coordinate-free, global and unique description of the attitude of the rigid body. Define the matrix 
composed of all $k$ measured vectors expressed in the body-fixed frame as column vectors,
\begin{align}
U^m&= [u_1^m\ u_2^m\ u_1^m\times u_2^m]  \mbox{ when }\, k=2,\, \mbox{ and }\nn\\ 
U^m&=[u_1^m\ u_2^m\ ... u_k^m]\in \mathbb{R}^{3\times k} \mbox{ when } k>2,\label{Umexp}
\end{align}
and the corresponding matrix of all these vectors expressed in the inertial frame as
\begin{align}
E&= [e_1\ e_2\ e_1\times e_2]  \mbox{ when }\, k=2,\, \mbox{ and }\nn\\ 
E&=[e_1\ e_2\ ... e_k]\in \mathbb{R}^{3\times k}  \mbox{ when } k>2.\label{Eexp}
\end{align}
Note that the matrix of the actual body vectors $u_j$ corresponding to the inertial vectors 
$e_j$, is given by
\begin{align}
U&= R \T E= [u_1\ u_2\ u_1\times u_2]  \mbox{ when }\, k=2,\, \mbox{ and }\nn\\ 
U&=R \T E=[u_1\ u_2\ ... u_k]\in \mathbb{R}^{3\times k} \mbox{ when } k>2.\nn
\end{align}
The direction vector measurements are given by
\begin{align}
u_j^m=R\T e_j+ \nu_j\, \mbox{ or }\, U^m = R\T E+ N, \label{DirMeasMod}
\end{align}
where $\nu_j\in\bR^3$ is an additive measurement noise vector and $N\in\bR^{3\times k}$ 
is the matrix with $\nu_j$ as its $j^{\mbox{th}}$ column vector.

The attitude kinematics for a rigid body is given by Poisson's equation:
\begin{align}
\dot{R}=R\Omega^\times, \label{kinematics}
\end{align}
where $\Omega\in\bR^3$ is the angular velocity vector and $(\cdot)^\times: \bR^3\to\so\subset
\bR^{3\times 3}$ is the skew-symmetric cross-product operator that gives a vector space 
isomorphism between $\bR^3$ and $\so$. The measurement model for angular velocity is 
\begin{align}
\Omega^m=\Omega+ w +\beta,\label{AngMeasMod}
\end{align}
where $w\in\bR^3$ is the measurement error in angular velocity and $\beta\in\bR^3$ is 
a vector of bias in angular velocity component measurements, which we consider to be a 
constant vector. 

\section{Attitude State and Bias Estimation Based on the Lagrange-d'Alembert 
Principle}\label{Sec4}
In order to obtain attitude state estimation schemes from continuous-time vector and angular 
velocity measurements, we apply the Lagrange-d'Alembert principle to an action functional of 
a Lagrangian of the state estimate errors, with a dissipation term linear in the angular velocity 
estimate error. This section presents an estimation scheme obtained using this approach. Let 
$\hat R\in\SO$ denote the estimated rotation matrix. According to \cite{Automatica}, the potential 
``energy'' function representing the attitude estimate error can be expressed as a generalized 
Wahba's cost function as
\begin{align}
\cU (\hat R,\U)=\Phi \Big(\frac{1}{2}\lan E-\hat R \U,(E-\hat R \U)W\ran \Big), \label{attindex}
\end{align}
where $\U$ is given by equation \eqref{Umexp}, $E$ is given by \eqref{Eexp}, and 
$W$ is the positive diagonal matrix of the weight factors for the measured 
directions. Note that $W$ may be generalized to be any positive 
definite matrix, not necessarily diagonal. Furthermore, $\Phi: [0,\infty)\mapsto[0,\infty)$ is a $C^2$ function that satisfies $\Phi(0)=0$ and 
$\Phi'(x)>0$ for all $x\in[0,\infty)$. Also $\Phi'(\cdot)\leq\alpha(\cdot)$ where 
$\alpha(\cdot)$ is a Class-$\mathcal{K}$ function. Let $\hat\Omega\in\bR^3$ 
and $\hat\beta\in\bR^3$ denote the estimated angular velocity and bias vectors, 
respectively. The ``energy" contained in the vector error between the estimated and the 
measured angular velocity is then given by
\be \cT (\hat\Omega,\Omega^m,\hat\beta)=\frac{m}{2}(\Omega^m-\hat\Omega-\hat\beta) \T 
(\Omega^m-\hat\Omega-\hat\beta), \label{angvelindex} \ee
where $m$ is a positive scalar. One can consider the Lagrangian composed of these 
``energy" quantities, as follows:
\begin{align} 
\cL (\hat R,\U,&\hat\Omega,\Omega^m,\hat\beta) = \cT(\hat\Omega,\Omega^m,\hat\beta)-
\cU (\hat R,\U) \nn\\
=&\frac{m}{2}(\Omega^m-\hat\Omega-\hat\beta) \T (\Omega^m-\hat\Omega-\hat\beta)\nn\\
&-\Phi\Big( \frac12\lan E-\hat{R}\U,(E-\hat{R}\U)W\ran\Big). \label{cLag}
\end{align}
If the estimation process is started at time $t_0$, then the action functional 
of the Lagrangian \eqref{cLag} over the time duration $[t_0,T]$ is expressed as
\begin{align}
\cS (\cL(\hat R&,\U,\hat\Omega,\Omega^m))= \int_{t_0}^T \big(\cT (\hat\Omega,
\Omega^m,\hat\beta)- \cU (\hat R,\U)\big)\di s \nn\\
=& \int_{t_0}^T \bigg\{ \frac{m}{2}(\Omega^m-\hat\Omega-\hat\beta) \T (\Omega^m-
\hat\Omega-\hat\beta)\nn\\
 &- \Phi\Big(\frac{1}{2}\lan E-\hat{R}\U,(E-\hat{R}\U)W\ran\Big) \bigg\} \di s.
\label{eq:J6}
\end{align}
Define the angular velocity measurement residual and the dissipation term:
\be \omega := \Omega^m-\hat\Omega- \hat\beta, \;\ \tau_D= D\omega, \label{angvelres} \ee
where $D\in\bR^{3\times 3}$ is positive definite.
Consider attitude state estimation in continuous time in the presence of measurement noise 
and initial state estimate errors. Applying the Lagrange-d'Alembert principle to the action 
functional $\cS (\cL(\hat R,\U,\hat\Omega,\Omega^m))$ given by \eqref{eq:J6}, in the 
presence of a dissipation term linear in $\omega$, leads to the 
following attitude and angular velocity filtering scheme. 
\begin{theorem} \label{filterN}
The filter equations for a rigid body with the attitude kinematics \eqref{kinematics} and with 
measurements of vectors and angular velocity in a body-fixed frame, are of the form
\begin{align}
\begin{cases}
&\dot{\hat{R}}=\hat{R}\hat{\Omega}^\times=\hat{R}(\Omega^m-\omega-\hat\beta)^\times,
\vspace{3mm}\\
&m\dot{\omega}= -m\hat{\Omega}\times \omega+\Phi'\big(\cU^0(\hat{R},\U)\big)S_L(\hat{R})-D\omega,\vspace{3mm}\\
&\hat\Omega=\Omega^m-\omega-\hat\beta,
\end{cases}
\label{eq:filterNoise}
\end{align}
where $D$ is a positive definite filter gain matrix, $\hat{R}(t_0)=\hat{R}_0$, $\omega(t_0)=\omega_0
=\Omega^m_0-\hat\Omega_0$, $S_L(\hat R)= \mrm{vex}\big(L\T \hat R - \hat R \T L\big)\in\bR^3$, $\mrm{vex}(\cdot): 
\so\to\bR^3$ is the inverse of the $(\cdot)^\times$ map, 
$L=EW(\U)\T$ and $W$ is chosen to satisfy the conditions in Lemma 2.1 of \cite{Automatica}.
\label{filter1}
\end{theorem}
{\em Proof}: In order to find an estimation scheme that filters the measurement noise in the 
estimated attitude, take the first variation of the action functional \eqref{eq:J6} with respect to 
$\hat R$ and $\hat\Omega$ and apply the Lagrange-d'Alembert principle with the dissipative 
term in \eqref{angvelres}. Consider the potential term $\cU^0(\hat R,\U)=\frac{1}{2}\lan E-\hat R\U,
(E-\hat R\U)W \ran$. 
Taking the first variation of this function with respect to $\hat{R}$ gives
\begin{align}
\delta\cU^0&=\lan -\delta\hat R\U,(E-\hat R\U)W \ran\nn\\
&=\frac{1}{2}\lan \Sigma^\times, \U WE\T\hat R-\hat R\T EW(\U)\T \ran,\nn\\
&=\frac{1}{2}\lan \Sigma^\times, L\T\hat R-\hat R\T L \ran=S\T_L(\hat R)\Sigma.
\end{align}
Now consider $\cU(\hat R,\U)=\Phi\big(\cU^0(\hat R,\U)\big)$. Then,
\begin{align}
\delta\cU=\Phi'\big(\cU^0(\hat R,\U)\big)\delta\cU^0=\Phi'\big(\cU^0(\hat R,\U)\big)S\T_L(\hat R)\Sigma.
\end{align}
Taking the first variation of the kinetic energy-like term \eqref{angvelindex} with respect to $\hat\Omega$ 
yields
\begin{align}
\delta\cT&=-m(\Omega^m-\hat\Omega-\hat\beta)\T\delta\hat\Omega \nn \\
&=-m(\Omega^m-\hat\Omega
-\hat\beta)\T(\dot\Sigma+\hat\Omega\times\Sigma)\nn\\
&=-m\omega\T(\dot\Sigma+\hat\Omega\times\Sigma),
\end{align}
where $\omega$ is as given by \eqref{angvelres}. 
Applying the Lagrange-d'Alembert principle and integrating by parts leads to
\begin{align}
&~~~~~\delta\cS+\int_{t_0}^T\tau_D\T\Sigma\di t=0\nn\\
&\Rightarrow -m\omega\T\Sigma\big|_{t_0}^T+\int_{t_0}^T m\dot{\omega}\T\Sigma\di t\\
&=\int_{t_0}^T\Big\{m\omega\T\hat\Omega^\times+\Phi'\big(\cU^0(\hat R,\U)\big)S\T_L(\hat R)-\tau_D\T\Big\}\Sigma\di t,\nn
\end{align}
where the first term in the left hand side vanishes, since $\Sigma(t_0)=\Sigma(T)=0$. After substituting 
$\tau_D=D\omega$, one obtains the second equation in \eqref{eq:filterNoise}.
\hfill\ensuremath{\square}

\section{Stability and Convergence of Variational Attitude Estimator}\label{Sec5}
The variational attitude estimator given by Theorem \ref{filterN} can be used in the presence  
of bias in the angular velocity measurements given by the measurement model 
\eqref{AngMeasMod}. The following analysis gives the stability and convergence properties of 
this estimator for the case that $\beta$ in \eqref{AngMeasMod} is constant.

\subsection{Stability of Variational Attitude Estimator} 
Prior to analyzing the stability of this attitude estimator, it is useful and instructive to interpret 
the energy-like terms used to define the Lagrangian in equation \eqref{cLag} in terms of state 
estimation errors. The following result shows that the measurement residuals, and therefore 
these energy-like terms, can be expressed in terms of state estimation errors. 
\begin{proposition}
Define the state estimation errors
\begin{align} 
&Q= R \hat R\T \, \mbox{ and }\, \omega= \Omega-\hat\Omega-\tilde\beta, \label{esterrs} \\
&\mbox{where }\, \tilde\beta = \beta -\hat\beta. \label{biaserr} 
\end{align}
In the absence of measurement noise, the energy-like terms \eqref{attindex} and 
\eqref{angvelindex} can be expressed in terms of these state estimation errors as follows:
\begin{align}
&\cU (Q) =\Phi \Big(\lan I-Q, K\ran \Big)\, \mbox{ where }\, K=EWE\T, \label{attindperf} \\
&\mbox{and } \cT (\omega)= \frac{m}{2} \omega\T\omega. \label{angvindperf} 
\end{align} 
\end{proposition}
{\em Proof}: The proof of this statement is obtained by first substituting $N=0$ and $w=0$ in 
\eqref{DirMeasMod} and \eqref{AngMeasMod}, respectively. The resulting expressions 
for $\U$ and $\Omega^m$ are then substituted back into \eqref{attindex} and 
\eqref{angvelindex}, respectively. Note that the same variable $\omega$ is used to represent 
the angular velocity estimation error for both cases: with and without measurement noise. 
Expression \eqref{attindperf} is also derived in \cite{Automatica}. 
\hfill\ensuremath{\square}
 
The stability of this estimator, for the case of constant rate gyro bias vector $\beta$, is given by  
the following result.
\begin{theorem}\label{stabproof}
Let $\beta$ in equation \eqref{AngMeasMod} be a constant vector. Then the variational attitude 
estimator given by equations \eqref{eq:filterNoise}, in addition to the following equation for 
update of the bias estimate:
\be \dot{\hat\beta}=  \Phi' \big(\cU^0(\hat{R},\U)\big) P^{-1} S_L(\hat{R}), \label{biasest} \ee
is Lyapunov stable for $P\in\bR^{3\times 3}$ positive definite.
\end{theorem}
{\em Proof}: To show Lyapunov stability, the following Lyapunov function is used:
\begin{align} 
V (\U,&\Omega^m,\hat{R},\hat\Omega,\hat\beta) = \frac{m}{2} (\Omega^m-\hat\Omega-\hat\beta)\T
(\Omega^m-\hat\Omega-\hat\beta) \nn\\
&+ \Phi \big(\cU^0(\hat{R},\U)\big) + \frac12 (\beta-\hat\beta)\T P (\beta-\hat\beta). 
\label{Lyapf} 
\end{align}  
Now consider the case that there is no measurement noise, i.e., $N=0$ and $w=0$ in 
equations \eqref{DirMeasMod} and \eqref{AngMeasMod}, respectively. In this case, 
the Lyapunov function \eqref{Lyapf} can be re-expressed in terms of the errors $\omega$, $Q$ and 
$\tilde\beta$ defined by equations \eqref{esterrs}-\eqref{biaserr}, as follows:
\be V(Q,\omega,\tilde\beta)= \frac{m}{2}\omega\T\omega + \Phi\big(\lan I-Q, K\ran\big) +
\frac12\tilde\beta\T P\tilde\beta. \label{errLyap} \ee
The time derivative of the attitude estimation error, $Q\in\SO$, is obtained as:
\be \dot Q= R(\Omega-\hat\Omega)^\times\hat{R}\T= Q\big(\hat R (\omega-\tilde\beta)\big)^\times, 
\label{dotQ} \ee
after substituting for $\hat\Omega$ from the third of equations \eqref{eq:filterNoise} in the 
case of zero angular velocity measurement noise (when $\Omega^m= \Omega+ \beta$). The 
time derivative of the Lyapunov function expressed as in \eqref{errLyap} can now be obtained 
as follows: 
\begin{align}
\dot V(Q,&\omega,\tilde\beta)= \label{dotV} \\
&m \omega\T\dot\omega -\Phi'\big(\lan I-Q, K\ran\big) S_L\T (\hat R) (\omega-\tilde\beta) -\tilde\beta\T P\dot{\hat\beta}.\nn
\end{align}
After substituting equation \eqref{biasest} and the second of equations \eqref{eq:filterNoise} in the 
above expression, one can simplify the time derivative of this Lyapunov function along the 
dynamics of the estimator as
\be \dot V(Q,\omega,\tilde\beta)= -\omega\T D\omega \le 0. \label{Vdotsimp} \ee
This time derivative is negative semi-definite in the estimate errors $(Q,\omega,
\tilde\beta)\in\Ta\SO\times\bR^3$. This proves the result. \hfill\ensuremath{\square}

\subsection{Domain of Convergence of Variational Attitude Estimator}
The domain of convergence of this estimator is given by the following result. 
\begin{theorem}\label{convproof}
The variational attitude estimator in the case of biased velocity measurements, given by eqs. 
\eqref{eq:filterNoise} and \eqref{biasest}, converges asymptotically to $(Q,\omega,\tilde\beta)=(I,0,
0)\in\Ta\SO\times\bR^3$ with an almost global domain of convergence.
\end{theorem}

The proof of this result is similar to the proof of the domain of convergence of the variational 
attitude estimator for the bias-free case in~\cite{Automatica}. The additional estimate error 
state $\tilde\beta$ converges to zero asymptotically for almost all initial $(Q,\omega)$ except 
those that lie on a set whose complement is dense and open in $\Ta\SO\simeq\SO\times\bR^3$. 

\section{Discrete-Time Estimator}\label{Sec6}
The ``energy" in the measurement residual for attitude is discretized as:
\begin{align}
\cU (\hat R_i,\U_i)&=\Phi\Big(\cU^0(\hat R_i,\U_i)\Big)\label{dattindex}\\
&=\Phi \Big(\frac{1}{2}\lan E_i-\hat R_i \U_i,(E_i-\hat R_i \U_i)W_i\ran \Big),\nn
\end{align}
where $W_i$ is a positive diagonal matrix of weight factors for the measured directions at time 
$t_i$, and $\Phi: [0,\infty)\mapsto[0,\infty)$ is a $C^2$ function that satisfies $\Phi(0)=0$ and 
$\Phi'(x)>0$ for all $x\in[0,\infty)$. Furthermore, $\Phi'(\cdot)\leq\alpha(\cdot)$ where 
$\alpha(\cdot)$ is a Class-$\mathcal{K}$ function.
The ``energy" in the angular velocity measurement residual is discretized as
\be \cT (\hat\Omega_i,\Omega^m_i)=\frac{m}{2}(\Omega^m_i-\hat\Omega_i-\hat\beta_i) \T
(\Omega^m_i-\hat\Omega_i-\hat\beta_i), \label{dangvelindex} \ee 
where $m$ is a positive scalar. 

Similar to the continuous-time attitude estimator in \cite{Automatica}, one can express 
these ``energy" terms for the case that perfect measurements (with no measurement noise) are 
available. In this case, these ``energy" terms can be expressed in terms of the state estimate 
errors $Q_i= R_i \hat R_i\T$ and $\omega_i= \Omega_i-\hat\Omega_i-\hat\beta_i$: 
\begin{align}
\begin{split}
& \cU (Q_i)= \Phi \Big(\frac{1}{2}\lan E_i - Q_i\T E_i,(E_i - Q_i\T E_i)W_i\ran \Big)=\\ 
&\Phi \big( \lan I-Q_i, K_i\ran \big)\, \mbox{ where }\, K_i= E_i W_i E_i\T, \\
& \mbox{and }\,\cT (\omega_i)= \frac{m}{2} \omega_i\T\omega_i\, \mbox{ where }
m>0. 
\end{split} \label{discUandT} 
\end{align}
The weights in $W_i$ can be chosen such that $K_i$ is always positive definite with distinct 
(perhaps constant) eigenvalues, as in the continuous-time estimator of \cite{Automatica}. Using 
these ``energy" terms in the state estimate errors, the discrete-time Lagrangian is expressed as:
\begin{align}
\cL (Q_i,\omega_i)&= \cT (\omega_i)- \cU ( Q_i) \nn\\
= &\frac{m}{2} \omega_i\T\omega_i- \Phi\big( \lan I-Q_i, K_i\ran \big).
\label{discLag}
\end{align}
The following statement gives a first-order discretization, in the form of a Lie 
group variational integrator, for the continuous-time estimator of Theorem \ref{filterN}.
\begin{proposition} \label{discfilter}
Let discrete-time measurements for two or more inertial vectors along with angular velocity be 
available at a sampling period of $h$. Further, let the weight 
matrix $W_i$ for the set of vector measurements $E_i$ be chosen such that $K_i=E_i 
W_i E_i\T$ satisfies Lemma 2.1 in \cite{Automatica}. A discrete-time estimator obtained by 
applying the discrete Lagrange-d'Alembert principle to the Lagrangian \eqref{discLag} is:
\begin{align}
&\hat R_{i+1}=\hat{R_i}\exp\big(h(\Omega_i^m-\omega_i
-\hat\beta_i)^\times\big),\label{1stDisFil_Rhat}\\
&\hat\beta_{i+1}= \hat\beta_i+ h\Phi' \big(\cU^0(\hat{R}_i,\U_i)\big) P^{-1} S_{L_i}(\hat{R}_i), 
\label{1stDisFil_betahat}\\
&\hat\Omega_i=\Omega_i^m-\omega_i-\hat\beta_i, \label{1stDisFil_Omegahat}\\
&m\omega_{i+1}=\exp(-h \hat\Omega_{i+1}^\times)\Big\{(m I_{3\times3}-hD)
\omega_i\label{1stDisFil_omega}\\
&~~~~~~~~~~~~~~~~~~~~~~~+h\Phi'\big(\cU^0(\hat R_{i+1},\U_{i+1})\big)S_{L_{i+1}}(\hat R_{i+1})\Big\},\nn
\end{align}
where $S_{L_i}(\hat R_i)=\mrm{vex}(L_i\T\hat R_i-\hat R_i\T L_i)\in\bR^3$, $L_i=E_i W_i(\U_i)\T\in
\mathbb{R}^{3\times3}$, $\cU^0(\hat R_i,\U_i)$ is defined in \eqref{dattindex} and $(\hat R_0,\hat\Omega_0)\in\SO\times\bR^3$ are initial estimated states.
\end{proposition}
The proof is very similar to the proof of the discrete-time variational attitude estimator presented in 
\cite{Automatica}. Note that the discrete-time total energy corresponding to the discrete Lagrangian 
\eqref{discLag} is dissipated with time, as with the continuous time estimator.

\section{Numerical Simulation}\label{Sec7}
This section presents numerical simulation results of the discrete estimator 
presented in Section \ref{Sec6}, in the presence of constant bias in angular velocity measurements. 
In order to validate the performance of this estimator, ``true" rigid body attitude states are generated 
using the rotational kinematics and dynamics equations. The rigid body moment of inertia is selected as 
$J_v=\diag([2.56\;\;3.01\;\;2.98]\T)$ kg.m$^2$. Moreover, a sinusoidal external torque is applied to this 
body, expressed in body fixed frame as
\begin{align}
\varphi(t)=[0\;\;\; 0.028\sin(2.7t-\frac{\pi}{7})\;\;\; 0]\T\mbox{ N.m}.
\end{align}
The true initial attitude and angular velocity are given by,
\begin{align}
\begin{split}
R_0=&\expm_{\SO}\bigg(\Big(\frac{\pi}{4}\times[\frac{3}{7}\;\;\;\; \frac{6}{7}\;\;\;\; \frac{2}{7}]\T\Big)^\times\bigg)\\
\mbox{and }& \Omega_0=\frac{\pi}{60}\times[-2.1\;\;\;\; 1.2\;\; -1.1]\T\mbox{ rad/s}.
\end{split}
\end{align}
A set of at least two inertial sensors and three gyros perpendicular to each other are assumed to be 
onboard the rigid body. The true states generated from the kinematics and dynamics of this rigid body 
are also used to generate the observed directions in the body fixed frame.
We assume that there are at most nine inertially known directions which are being measured 
by the sensors fixed to the rigid body at a constant sample rate. 
Bounded zero mean noise is added to the true direction vectors to generate each measured 
direction. A summation of 
three sinusoidal matrix functions is added to the matrix $U$, to generate a measured $\U$ 
with measurement noise. The frequency of the noises are 1, 10 and 100 Hz, with different phases 
and different amplitudes, which are up to $2.4^\circ$ based on coarse attitude sensors like sun 
sensors and magnetometers. Similarly, two sinusoidal noises of 10 Hz and 200 Hz frequencies 
are added to $\Omega$ to form the measured $\Omega^m$. These signals also have different 
phases and their magnitude is up to $0.97^\circ/s$, which corresponds to a coarse 
rate gyro. Besides, the gyro readings are assumed to contain a constant bias in three directions, as follows:
\begin{align}
\beta=[-0.01\;\;\;-0.005\;\;\;\;\;0.02]\T\mbox{ rad/s}.
\end{align}
The estimator is simulated over a time interval of $T=40$s, with a time stepsize of $h=0.01$s. The scalar
inertia-like gain is $m=5$ and the dissipation matrix is selected as: 
\begin{align}
D=\diag\big([17.4\;\;\; 18.85\;\;\; 20.3]\T\big).
\end{align}
As in \cite{Automatica}, $\Phi(x)=x$. The weight matrix $W$ is also calculated using the conditions in 
\cite{Automatica}. The positive definite matrix for bias gain is selected as $P=2\times10^3I$. The initial 
estimated states and bias are set to:
\begin{align}
\begin{split}
\hat R_0&=\expm_{\SO}\bigg(\Big(\frac{\pi}{2.5}\times[\frac{3}{7}\;\;\;\; \frac{6}{7}\;\;\;\; \frac{2}{7}]\T\Big)^\times\bigg),\\
\hat\Omega_0&=[-0.26\;\;\;\;\; 0.1725\;\;\; -0.2446]\T\mbox{ rad/s},\\
&\mbox{and }\hat\beta_0=[0\;\;\;-0.01\;\;\;\;\;0.01]\T\mbox{ rad/s}.
\end{split}
\end{align}

In order to integrate the implicit set of equations in \eqref{1stDisFil_Rhat}-\eqref{1stDisFil_omega} 
numerically, the first two equations are solved at each sampling step. Using \eqref{1stDisFil_Omegahat}, 
$\hat\Omega_{i+1}$ in \eqref{1stDisFil_omega} is written in terms of $\omega_{i+1}$ next. The resulting 
implicit equation is solved with respect to $\omega_{i+1}$ iteratively to a set tolerance applying the 
Newton-Raphson method. The root of this nonlinear equation along with $\hat R_{i+1}$ and 
$\hat\beta_{i+1}$ are used for the next sampling time instant. This process is 
repeated till the end of the simulated duration.\par
\begin{figure}
\begin{center}
\includegraphics[height=2.4in]{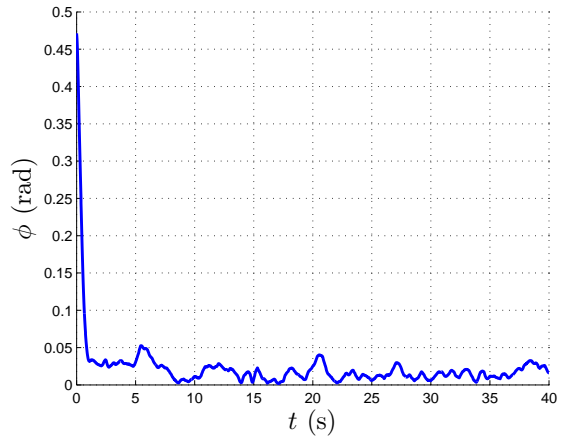}
\caption{Principal angle of the attitude estimate error}  
\label{Fig1}                                 
\end{center}                                 
\end{figure}
\begin{figure}
\begin{center}
\includegraphics[height=2.4in]{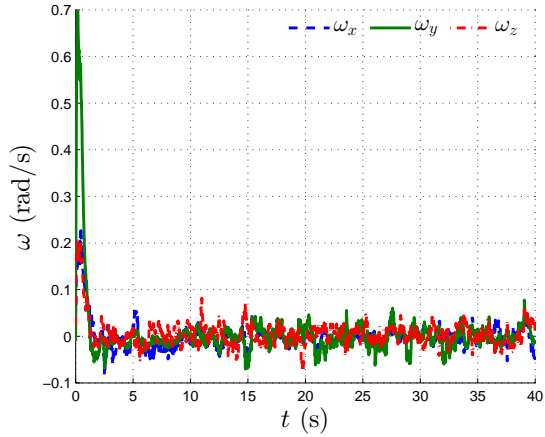}
\caption{Angular velocity estimate error}  
\label{Fig2}                                 
\end{center}                                 
\end{figure}
\begin{figure}
\begin{center}
\includegraphics[height=2.4in]{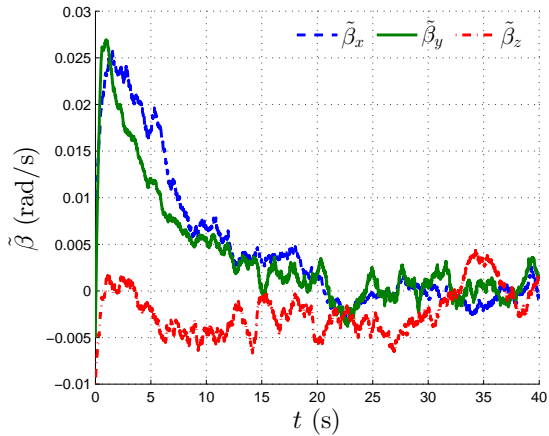}
\caption{Bias estimate error}  
\label{Fig3}                                 
\end{center}                                 
\end{figure}

Results from this numerical simulation are shown here. The principal angle corresponding 
to the rigid body's attitude estimation error is depicted in Fig. \ref{Fig1}, and estimation errors in the 
angular velocity components are shown in Fig. \ref{Fig2}. Finally, Fig. \ref{Fig3} portrays estimate errors in 
bias components. Estimation errors are seen to converge to a neighborhood of $(Q,\omega,\tilde\beta)
=(I,0,0)$, where the size of this neighborhood depends on the bounds of the measurement noise.

\section{Conclusion}\label{Sec8}
The framework of variational attitude estimation is generalized to include bias in angular 
velocity measurements and estimate a constant bias vector. The continuous-time state 
estimator is obtained by applying the Lagrange-d'Alembert principle from variational mechanics 
to a Lagrangian consisting of the energies in the measurement residuals, along with a 
dissipation term linear in angular velocity measurement residual. The update law for the 
bias estimate ensures that the total energy content in the state and bias estimation errors 
is dissipated as in a dissipative mechanical system. The resulting generalization of the 
variational attitude estimator is almost globally asymptotically stable, like the variational attitude 
estimator for the bias-free case reported in~\cite{Automatica}. A discretization of this estimator 
is obtained in the form of an implicit first order Lie group variational integrator, by applying the 
discrete Lagrange-d'Alembert principle to the discrete Lagrangian with the dissipation term 
linear in the angular velocity estimation error. This discretization preserves the stability of 
the continuous estimation scheme. Using a realistic set of data for rigid body 
rotational motion, numerical simulations show that the estimated states and estimated bias 
converge to a bounded neighborhood of the true states and true bias when the measurement 
noise is bounded. Future planned extensions of this work are to develop an explicit 
discrete-time implementation of this attitude estimator, 
and implement it in real-time with optical and inertial sensors.  


\bibliographystyle{IEEEtran}        
\bibliography{alias,VAE_refs}

\begin{thebibliography}{10}
\providecommand{\url}[1]{#1}
\csname url@samestyle\endcsname
\providecommand{\newblock}{\relax}
\providecommand{\bibinfo}[2]{#2}
\providecommand{\BIBentrySTDinterwordspacing}{\spaceskip=0pt\relax}
\providecommand{\BIBentryALTinterwordstretchfactor}{4}
\providecommand{\BIBentryALTinterwordspacing}{\spaceskip=\fontdimen2\font plus
\BIBentryALTinterwordstretchfactor\fontdimen3\font minus
  \fontdimen4\font\relax}
\providecommand{\BIBforeignlanguage}[2]{{%
\expandafter\ifx\csname l@#1\endcsname\relax
\typeout{** WARNING: IEEEtran.bst: No hyphenation pattern has been}%
\typeout{** loaded for the language `#1'. Using the pattern for}%
\typeout{** the default language instead.}%
\else
\language=\csname l@#1\endcsname
\fi
#2}}
\providecommand{\BIBdecl}{\relax}
\BIBdecl

\bibitem{Automatica}
M.~Izadi and A.~Sanyal, ``Rigid body attitude estimation based on the
  {L}agrange-d'{A}lembert principle,'' \emph{Automatica}, vol.~50, no.~10, pp.
  2570 -- 2577, 2014.

\bibitem{ICRA2015}
M.~Izadi, A.~Sanyal, E.~Samiei, and V.~Kumar, ``Comparison of an attitude
  estimator based on the {L}agrange-d'{A}lembert principle with some
  state-of-the-art filters,'' in \emph{2015 IEEE International Conference on
  Robotics and Automation, May 26 - 30, 2015}, Seattle, Washington, 2015.

\bibitem{leishman2014quadrotors}
R.~Leishman, J.~Macdonald, R.~Beard, and T.~McLain, ``Quadrotors and
  accelerometers: State estimation with an improved dynamic model,''
  \emph{Control Systems, IEEE}, vol.~34, no.~1, pp. 28--41, 2014.

\bibitem{bras2013nonlinear}
S.~Br{\'a}s, M.~Izadi, C.~Silvestre, A.~Sanyal, and P.~Oliveira, ``Nonlinear
  observer for {3D} rigid body motion,'' in \emph{Decision and Control (CDC),
  2013 IEEE 52nd Annual Conference on}.\hskip 1em plus 0.5em minus 0.4em\relax
  IEEE, 2013, pp. 2588--2593.

\bibitem{morgado2014embedded}
M.~Morgado, P.~Oliveira, C.~Silvestre, and J.~Vasconcelos, ``Embedded vehicle
  dynamics aiding for {USBL/INS} underwater navigation system,'' \emph{Control
  Systems Technology, IEEE Transactions on}, vol.~22, no.~1, pp. 322--330,
  2014.

\bibitem{gold}
H.~Goldstein, \emph{Classical Mechanics}, 2nd~ed.\hskip 1em plus 0.5em minus
  0.4em\relax Reading, MA: Addison-Wesley, 1980.

\bibitem{green}
D.~Greenwood, \emph{Classical Dynamics}.\hskip 1em plus 0.5em minus 0.4em\relax
  Englewood Cliffs, NJ: Prentice Hall, 1987.

\bibitem{marswest}
J.~Marsden and M.~West, ``Discrete mechanics and variational integrators,''
  \emph{Acta Numerica}, vol.~10, pp. 357--514, 2001.

\bibitem{haluwa}
E.~Hairer, C.~Lubich, and G.~Wanner, \emph{Geometric Numerical
  Integration}.\hskip 1em plus 0.5em minus 0.4em\relax New York: Springer
  Verlag, 2002.

\bibitem{black64}
H.~Black, ``A passive system for determining the attitude of a satellite,''
  \emph{AIAA Journal}, vol.~2, no.~7, pp. 1350--1351, 1964.

\bibitem{jo:wahba}
G.~Wahba, ``A least squares estimate of satellite attitude, {P}roblem 65-1,''
  \emph{SIAM Review}, vol.~7, no.~5, p. 409, 1965.

\bibitem{jo:solwahba}
J.~Farrell, J.~Stuelpnagel, R.~Wessner, J.~Velman, and J.~Brock, ``A least
  squares estimate of satellite attitude, {S}olution 65-1,'' \emph{SIAM
  Review}, vol.~8, no.~3, pp. 384--386, 1966.

\bibitem{markley1988attitude}
F.~Markley, ``Attitude determination using vector observations and the singular
  value decomposition,'' \emph{The Journal of the Astronautical Sciences},
  vol.~36, no.~3, pp. 245--258, 1988.

\bibitem{sanyal2006optimal}
A.~Sanyal, ``Optimal attitude estimation and filtering without using local
  coordinates, {P}art 1: Uncontrolled and deterministic attitude dynamics,'' in
  \emph{American Control Conference, 2006}, Minneapolis, MN, 2006, pp.
  5734--5739.

\bibitem{bhat}
S.~P. Bhat and D.~S. Bernstein, ``A topological obstruction to continuous
  global stabilization of rotational motion and the unwinding phenomenon,''
  \emph{Systems \& Control Letters}, vol.~39, no.~1, pp. 63--70, 2000.

\bibitem{CSMpaper}
N.~A. Chaturvedi, A.~K. Sanyal, and N.~H. McClamroch, ``Rigid body attitude
  control---{U}sing rotation matrices for continuous, singularity-free control
  laws,'' \emph{{IEEE} Control Systems Magazine}, vol.~31, no.~3, pp. 30--51,
  2011.

\bibitem{jgcd12}
A.~Sanyal and N.~Nordkvist, ``Attitude state estimation with multi-rate
  measurements for almost global attitude feedback tracking,'' \emph{{AIAA}
  Journal of Guidance, Control, and Dynamics}, vol.~35, no.~3, pp. 868--880,
  2012.

\bibitem{silvest08}
J.~F. Vasconcelos, C.~Silvestre, and P.~Oliveira, ``A nonlinear {GPS/IMU} based
  observer for rigid body attitude and position estimation,'' in \emph{{IEEE}
  Conf. on Decision and Control}, Cancun, Mexico, Dec. 2008, pp. 1255--1260.

\bibitem{Lageman}
C.~Lageman, J.~Trumpf, and R.~Mahony, ``Gradient-like observers for invariant
  dynamics on a {L}ie group,'' \emph{IEEE Transaction on Automatic Control},
  vol.~55, pp. 367 -- 377, 2010.

\bibitem{markSO3}
F.~Markley, ``Attitude filtering on {SO}(3),'' \emph{The Journal of the
  Astronautical Sciences}, vol.~54, no.~4, pp. 391--413, 2006.

\bibitem{mahapf08}
R.~Mahony, T.~Hamel, and J.-M. Pfimlin, ``Complementary filters on the special
  orthogonal group,'' \emph{IEEE Transactions on Automatic Control}, vol.~53,
  no.~5, pp. 1203--1217, 2008.

\bibitem{bonmaro09}
S.~Bonnabel, P.~Martin, and P.~Rouchon, ``Nonlinear symmetry-preserving
  observers on {L}ie groups,'' \emph{IEEE Transactions on Automatic Control},
  vol.~54, no.~7, pp. 1709--1713, 2009.

\bibitem{Vas1}
J.~Vasconcelos, R.~Cunha, C.~Silvestre, and P.~Oliveira, ``A nonlinear position
  and attitude observer on {SE}(3) using landmark measurements,'' \emph{Systems
  \& Control Letters}, vol.~59, pp. 155--166, 2010.

\bibitem{Mortensen}
R.~Mortensen, ``Maximum-likelihood recursive nonlinear filtering,''
  \emph{Journal of Optimization Theory and Applications}, vol.~2, no.~6, pp.
  386--394, 1968.

\bibitem{aguhes06}
A.~Aguiar and J.~Hespanha, ``Minimum-energy state estimation for systems with
  perspective outputs,'' \emph{IEEE Transactions on Automatic Control},
  vol.~51, no.~2, pp. 226--241, 2006.

\bibitem{ZamPhD}
M.~Zamani, ``Deterministic attitude and pose filtering, an embedded {L}ie
  groups approach,'' Ph.D. dissertation, Australian National University,
  Canberra, Australia, Mar. 2013.

\bibitem{kirk}
D.~Kirk, \emph{Optimal Control Theory: An Introduction}.\hskip 1em plus 0.5em
  minus 0.4em\relax Englewood Cliffs, NJ: Prentice Hall, 1970.

\end{thebibliography}

\end{document}